\documentclass[12pt,a4paper]{article}
\usepackage{amssymb}
\usepackage{amsmath}
\newtheorem{theoreme}{Theorem}[section]
\newtheorem{prop}[theoreme]{Proposition}

\newcommand{\card}{{\rm card}}
\newcommand{\R}{\mathbb{R}}
\newcommand{\Q}{\mathbb{Q}}

\newcommand{\Z}{\mathbb{Z}}
\newcommand{\N}{\mathbb{N}}

\newcommand{\points}{{\cal M}_{\sigma}}

\begin{document}

\title{Diffraction and Palm measure of point processes}

\author{Jean-Baptiste Gou\'er\'e\footnote{\textit{Postal address}: Universit\'e Claude 
Bernard Lyon 1,
LaPCS, B\^at. B, Domaine de Gerland, 50, avenue Tony-Garnier F-69366 Lyon 
Cedex 07, France.
\textit{E-mail}: jbgouere@univ-lyon1.fr }}  

\date{}

\maketitle

\begin{abstract}
Using the Palm measure notion, we prove the existence of the diffraction measure of all stationary and ergodic
point processes. We get precise expressions of those measures in the case of specific processes : stochastic
subsets of $\Z^d$, sets obtained by the ``cut-and-project'' method.
\end{abstract}

$ $\\
From a physical point of view, the diffraction of X rays by a material is a way of studying its microscopic structure.
We remind the mathematical formalism of diffraction theory (Hof \cite{Hof-95}).
Let $\chi$ be a a locally finite subset of $\R^d$ (physically, each point of $\chi$ is the center of an atom).
The autocorrelation $\gamma_R$ is defined by
\begin{equation} \label{def-ac-points}
\gamma_R:=\frac{1}{|B_R|} \sum_{x,y \in \chi \cap B_R} \delta_{y-x}
=\frac{1}{|B_R|} \; \left(\sum_{x \in \chi \cap B_R} \delta_x\right)*\left(\sum_{x \in \chi \cap B_R} \delta_{-x}\right)
\end{equation}
where $\delta_x$ is the Dirac measure at $x$ and $|B_R|$ is the canonical Lebesgue measure of the ball $B_R$ of radius $R>0$ centered at $0$.
The autocorrelation $\gamma_R$, also called the Patterson function,
represents the relative positions of the different points of $\chi \cap B_R$.
Physically, the Fourier transform of that measure corresponds with the diffraction pattern : 
$\widehat{\gamma_R}(t)$ is the luminous intensity diffracted in the direction $t$ by a material whose atom centers
are the points of $\chi \cap B_R$.

If it exists, the limit $\gamma$ in the vague topology as $R \rightarrow \infty$ of the measures $\gamma_R$
is called the autocorrelation of $\chi$.
That measure is tempered and its Fourier transform $\widehat{\gamma}$ is a positive measure called the
diffraction measure of $\chi$.

\section{Autocorrelation and Palm measure}

Let $\points(\R^d)$ be the set of all locally finite subsets of $\R^d$.
With each bounded borel set $A$ we associate the application
$$N_A \; : \; 
\begin{cases}
\points & \longrightarrow \N  \\
\phi & \longmapsto \card(A \cap \phi). 
\end{cases}
$$
We endow $\points$ with the $\sigma$-algebra $\cal{A}$ generated by those applications.
A point process \cite{Moller,Neveu-pp} is a measurable application $\chi$ from a probability space
$(\Omega,\cal{F},P)$ to $(\points,\cal{A})$.
If we equip $(\points,\cal{A})$ with the law of $\chi$
we can suppose that $(\Omega,\cal{F})=(\points,\cal{A})$ and that  $\chi$ if the identity application.

A point process is said to be integrable if the random variables $N_A$ are integrable for
every compact $A$. It's said to be square integrable if those variables are square integrable.

If $t \in \R^d$ let $T_t$ denote the translation of $\points$
$$T_t \; : \; 
\begin{cases}
\points & \longrightarrow \points  \\
\phi={\{x_n\}}_n & \longmapsto \phi-t={\{x_n-t\}}_n.
\end{cases}
$$
A point process is said to be stationary if its law is invariant under the action of the translations
$(T_t)_{t \in \R^d}$.
Let $\chi$ be a stationary and integrable point process.
The Palm measure of $\chi$ is a measure $\widetilde{P}$ on $(\points,\cal{A})$ defined by  \cite{Moller}
$$\widetilde{P}(F) := \frac{1}{|B|}E\sum_{x \in \phi \cap B} 1_F(\phi-x), \quad F \in \cal{A}$$
where $B$ is a fixed borel subset of $\R^d$ whose Lebesgue measure $|B|$ is finite and strictly positive.
That definition does not depend on the choice of $B$.

With each measure $m$ on $(\points,{\cal A})$ we associate the intensity measure $I(m)$ on  $(\R^d,{\cal B}(\R^d))$ 
defined by
$$I(m)(A)=\int \card(\phi \cap A) \; dm(\phi), \quad  A \in {\cal B}(\R^d)$$
where ${\cal B}(\R^d)$ is the set of borel subsets of $\R^d$.

\begin{theoreme} \label{ac-palm}
Let $\chi$ be a stationary ergodic and square integrable point process on
$\mathbb{R}^d$. Let $\widetilde{P}$ be its Palm measure and $I(\widetilde{P})$ the intensity of that measure.
Then $I(\widetilde{P})$ is locally finite and
for each bounded borel set $C$ of $\R^d$ we have
$$\lim_{R \rightarrow \infty} \gamma_R(C) = I(\widetilde{P})(C)  \;\;\; a.e.$$
In particular, $\gamma_R$ converge a.e. to $I(\widetilde{P})$ in the vague topology.
\end{theoreme}
$ $\\
{\bf Example : } Let $\chi$ be a Poisson point process on $\R^d$ whose intensity measure is the canonical 
Lebesgue measure on $\R^d$ \cite{Moller}.
The conditions of the theorem \ref{ac-palm} are classicaly fulfilled and we know
\cite {Moller} that the Palm measure of $\chi$ is the law of the process $\chi_0:=\chi \cup \{0\}$.
We immediately deduce that the process $\chi$ admits a.e. an autocorrelation measure $\gamma$ which verifies
$\gamma=\widehat{\gamma}=\delta_0+dx$.

\section{Stochastic subsets of $\Z^d$} \label{section-reseau-aleat}

\subsection{{\bf Generalities}}

Let $X=(X_k)_{k \in \Z^d}$ be a stationary and ergodic process with values in
$\{0,1\}$. Let $\mu$ be its spectral measure.
We consider the process $\chi = \{k \in \Z^d \; : \; X_k=1\}$.

\begin{theoreme} \label{th-01} 
\begin{enumerate}
\item The process $\chi$ admits a.e. an autocorrelation $\gamma$. We have
\begin{equation} \label{ac-01}
\gamma := \sum_{k \in \Z^d} EX_0X_k\;\delta_k
\end{equation}
\item The diffraction measure is $\widehat{\gamma}=\mu_p$
where $\mu_p$ is the $\Z^d-$periodic measure associated with $\mu$ by
\begin{equation} \label{periodiser}
\mu_p(A)=\sum_{k \in \Z^d} \mu(A-k), \quad A \in {\cal B}(\R^d).
\end{equation}
\end{enumerate}
\end{theoreme}
{\bf Remark : }
Let $U$ be a uniform random variable on $[0,1[^d$. We assume that $U$ and $X$ are independant.
We build a ``$\R^d$-stationary'' version of $\chi$ writing
$\varphi=s(U,X)$ where $s$ if the function from $[0,1[^d\times \{0,1\}^{\Z^d}$ to $\points(\R^d)$ defined by 
\begin{equation} \label{def-f-01}
s(u,x)=u+\{k \in \Z^d \; : \; x_k=1\}.
\end{equation}
The results of the theorem \ref{th-01} are obviously verified by the process $\varphi$.

\subsection{An example of a singular diffraction measure} \label{ex-construction-01}

Let $f$ be the function defined on $[-1,1]$ by
\begin{equation} \label{f-signe}
f(x):=\frac{1}{2\pi}\arcsin(x)+\frac{1}{4}=\sum_{k \in \N} a_k x^k,
\end{equation}
where all $a_k$ are positive.
If $\nu$ is a probability measure on $\R^d$ we define
the measure $f(\nu)$ by $f(\nu)=\sum_{k \in \N} a_k \nu^{*k}$.

\begin{prop} \label{diff-per}
Let $\nu$ be a symetric and continue probability measure on $\R$. Write $\mu=f(\nu)$
where $f$ is defined by (\ref{f-signe}). There exists a point process which admits
a.e. the diffraction measure $\mu_p$ 
($\mu_p$ is associated with $\mu$ by (\ref{periodiser})). 
\end{prop}
{\bf Example : }
Let $\nu$ be the law of the random serie $\sum_{n \geq 1} \frac{1}{n^n} X_n$ where ${(X_n)}_{n \in \N}$ 
is a sequence of i.i.d. random variables with $P(X_i=1)=P(X_i=-1)=1/2$.
We get a point process which admits a diffraction measure $\mu$ which can be decomposed into $\mu=\frac{1}{4}\sum_{k \in \Z} \delta_k + m$
where $m=[(f-1/4)(\nu)]_p$ is a purely singular $1$-periodic measure.

\section{Sets built by the ``cut-and-project'' method} \label{section-cut-proj}

\subsection{Principle}

Let $E$ and $F$ be two linear subspaces of $\R^d$ such that $\R^d=E \oplus F$.
If $x \in \R^d$ we write $x=x_E+x_F=p_E(x)+p_F(x)$ with obvious notations.
Let  $x_E^{\perp}$ and $x_F^{\perp}$ be the orthogonal projections of $x$ on $E$ and $F$.
Let W be a fixed bounded borel subset of $F$.
Let $\phi$ be a stationary and square integrable point process of $\R^d$.
We assume that $\phi$ is ergodic under the action of the translations ${(T_t)}_{t \in E}$.
We define a new point process $\chi$ on $E$ by
$\chi=\pi(\phi)$ where $\pi$ is the application from  $\points(\R^d)$ to $\points(E)$ defined by
\begin{equation} \label{projection}
\pi(\phi)=p_E\left(\phi\cap (E \times W)\right).
\end{equation}
We make the following assumption :
\begin{equation} \label{hypothesealacon} 
\hbox{Two different points of }\phi \cap (E \times W)\hbox{ have a different image by }p_E\hbox{ (a.e.).}
\end{equation}

\begin{theoreme} \label{th-cut-proj}
Let $\lambda$ denote the restriction to $W$ of the canonical Lebesgue measure on $F$.
Let $\widetilde{Q}$ be the Palm measure of $\phi$.
Let $C_E$ be a borel subset of $E$ such that ${|C_E|}_E=1$ where ${|\cdot|}_E$ is the canonical Lebesgue measure on $E$.
We define $C_F$ similarly.
Let $\alpha$ denote the canonical Lebesgue measure (in $\R^d$) of $C_E+C_F$.
Then
\begin{enumerate}
\item The process $\chi$ is stationary ergodic and square integrable.
\item The Palm measure of $\chi$ is the image of the measure $\alpha \lambda \otimes \widetilde{Q}$ by the application
\begin{equation} \label{cut-proj-constr-palm}
f \; : \; 
\begin{cases}
W \times \points(\R^d) & \longrightarrow \points(E)  \\
(w,\widetilde{\phi}) & \longmapsto p_E\left((w+\widetilde{\phi})\cap (E \times W)\right)=\pi(w+\widetilde{\phi}).
\end{cases}
\end{equation}
\item The intensity of the Palm measure is the measure $\gamma$ on $E$ defined by
\begin{equation} \label{cut-proj-ac}
\gamma(C) = \alpha \int_{\R^d} I(\widetilde{Q})(dx) \; \psi(x_F)1_C(x_E), \quad C \in {\cal B}(E),
\end{equation}
where $\psi$ is the function from $F$ to $\R$ defined by $\psi(x) = 1_W * 1_{-W} (x)$.
\item
The Fourier transform of $\gamma$ is
\begin{equation} \label{cut-proj-diff}
\widehat{\gamma}(C)=
\alpha^2 \int_{\R^d} \widehat{I(\widetilde{Q})}(dx) \; \widehat{\psi}(x_F^{\perp}) 1_C(x_E^{\perp}), \quad C \in {\cal B}(E),
\end{equation}
where $\widehat{I(\widetilde{Q})}$ is the Fourier transform of $I(\widetilde{Q})$.
\end{enumerate}
\end{theoreme}
{\bf Remark : } With the formalism of marked point processes
(which allows the ``presence'' of several points at the same place)
we can avoid the assumption (\ref{hypothesealacon}). That changes neiher the proof nor the result.\\
$ $\\
We immediately deduce, with the theorem \ref{ac-palm}, the following result :

\begin{theoreme} \label{cor-cut-proj} 
The process $\chi$ admits a.e. the autocorrelation measure described by (\ref{cut-proj-ac}) and the diffraction measure described
by (\ref{cut-proj-diff}). 
\end{theoreme}

\subsection{Application to stochastic subsets of $\Z^d$} \label{section-reseau-cut}

We apply the construction ``cut-and-project'' to the process  $\varphi = s(U,X)$ introduced in the remark of the
section \ref{section-reseau-aleat}.
The case where the $X_k$ are i.i.d. has been studied in \cite{Baake-Moody-entropy}.

In order to apply the theorem \ref{th-cut-proj}, we have to verify the ergodicity of $\varphi$ under the action of the translations
${(T_t)}_{t \in E}$. To that effect, we prove that

\begin{prop} \label{prop-erg-direction}
There exists a negligeable set $S \subset S^{d-1}$ such that, for all  $v \in S^{d-1}\setminus S$,
$\varphi$ is ergodic under the direction $\R v$. If $d=2$ then $S$ is countable. 
\end{prop} 
In some particular cases, we can make $S$ explicit :
\begin{prop} Let ${(A^i_n)}_{n \in \Z}$, $1 \leq i \leq d$ be $d$ independant stationary and ergodic processes
taking values in $\{0,1\}$. We define a process $X$ on $\Z^d$ by
 $X_{n_1,..,n_d}=A^1_{n_1}..A^d_{n_d}$. That process is stationary and ergodic. If the processes
$A^i$ are weak-mixing (which is the case of the processes built in \ref{ex-construction-01}),
then we can take for $S$
the set of all vectors of $S^{d-1}$ whose coordinates are dependant over $\Q$.
\end{prop}

The process $\varphi$ is ergodic under the action of the translations ${(T_t)}_{t \in E}$ as soon as $E$ contains a vector $v \in S^{d-1}\setminus S$.
We have therefore the following result :
 
\begin{theoreme} \label{th-erg-direction}
Let $E$ and $F$ be two linear subspaces of $\R^d$ such that $\R^d=E \oplus F$.
We assume that $E$ contains a vector of $S^{d-1}\setminus S$.
Let $\chi=\pi(\varphi)$ where $\pi$ is defined by  (\ref{projection}).
The process $\chi$ is then stationary ergodic and square integrable.
It admits a.e. (with the notations of theorems \ref{th-cut-proj} and \ref{th-01}) an autocorrelation measure
such that
$$\gamma(C)=\alpha \sum_{k \in \Z^d} \widehat{\mu}(k)\psi(k_F)1_C(k_E)$$
and a diffraction measure such that
$$\widehat{\gamma}(C)=\alpha^2 \int_{\R^d} \widehat{\psi}(x_F^{\perp})1_C(x_E^{\perp}) \mu_p(dx), \quad C \in {\cal B}(\R^d).$$
\end{theoreme}
$ $\\

If we consider $X\equiv 1$ then the point process $\varphi$ is (modulo a translation) 
the deterministic lattice $\Z^d$. 
We can take for $S$ the set of all vectors of $S^{d-1}$ whose coordinates are dependant over $\Q$
and we then get with the theorem  \ref{th-erg-direction}
the classical results about ``model-sets''
\cite{Hof-95,Meyer-quasicrystals,Moody-unif-distribution,Schlottmann-book-Baake-Moody,Solomyak-spectrum-delone-set}.

\end{document}